\documentclass[12pt]{amsart}
\usepackage{graphicx,amsmath,amsfonts,amssymb,bbm}
\usepackage{psfrag}
\usepackage{epsf}

\makeatletter\@addtoreset {equation}{section}\makeatother

\newtheorem{theo}{Theorem}
\newtheorem{lem}{Lemma}[section]

\newtheorem{rem}{Remark}[section]
\newenvironment{Proof}
{\begin{trivlist} \item[]{\bf Proof. }}%
{\hspace*{\fill}$\rule{.3\baselineskip}{.35\baselineskip}$\end{trivlist}}

\renewcommand{\geq}{\geqslant}
\renewcommand{\leq}{\leqslant}

\renewcommand{\phi}{\varphi}
\newcommand{\be}{\begin{eqnarray}}
\newcommand{\ee}{\end{eqnarray}}

\newcommand{\eps}{\varepsilon}

\begin{document}

\title[Normal form for transverse instability]{Normal form for transverse instability of the line soliton with a nearly critical speed of propagation}

\author{Dmitry Pelinovsky}
\address{Department of Mathematics, McMaster University, Hamilton, Ontario  L8S 4K1, Canada}
\email{dmpeli@math.mcmaster.ca}
\address{Department of Applied Mathematics,
Nizhny Novgorod State Technical University, 24 Minin street, 603950 Nizhny Novgorod, Russia}

\thanks{The results of this work were obtained with the financial support from
the state task of Russian Federation in the sphere of scientific activity (Task No. 5.5176.2017/8.9).}

\date{\today}

\begin{abstract}
In the context of the line solitons in the Zakharov--Kuznetsov (ZK) equation,
there exists a critical speed of propagation such that small transversely periodic perturbations are unstable
if the soliton speed is larger than the critical speed and orbitally stable if the soliton speed is
smaller than the critical speed. The normal form for transverse instability
of the line soliton with a nearly critical speed of propagation is derived by means
of symplectic projections and near-identity transformations. Justification of this normal form
is provided with the energy method. The normal form predicts a transformation of the unstable
line solitons with larger-than-critical speeds
to the orbitally stable transversely modulated solitary waves.
\end{abstract}

\maketitle

\section{Introduction}

Starting with the pioneer works \cite{KP,ZR},
it is well known that the line solitons are spectrally unstable with respect to the long transverse perturbations
in many nonlinear evolution equations such as the Kadometsev--Petviashvili (KP) and
nonlinear Schr\"{o}dinger (NLS) equations (see review in \cite{kivpel}).
The spectral instability persists to the short transverse perturbations of any period in the hyperbolic
version of the two-dimensional NLS equation \cite{Dec1,Dec2}, whereas it disappears for
short transverse perturbations in the elliptic version of the two-dimensional NLS equation
and in the KP-I equation \cite{RT0,RT1}. Alternatively, for a fixed period of the transverse perturbation,
the transverse instability occurs for the line solitons with larger-than-critical speeds of propagation
and disappears for those with smaller-than-critical speeds.

In the prototypical case of the KP-I equation, it was shown in \cite{RT2} that if the line solitons
are spectrally stable with respect to the periodic perturbations, they remain nonlinearly stable,
whereas if they are spectrally unstable, they remain nonlinearly unstable. The spectral
stability analysis is inconclusive for the line soliton with the critical speed of propagation
since the linearized operator has an additional zero eigenvalue
beyond the one induced by the translational symmetry of the KP-I equation. The presence of the additional zero
eigenvalue implies a bifurcation of the new travelling solutions which are spatially localized
along the longitudinal direction and are periodic along the transverse perturbations.
Such travelling solitary waves with periodic transverse modulations were discovered
for the KP-I equation by Zaitsev \cite{Zaitsev}.

Analytical solutions for the unstable eigenmode were derived for the KP-I equation
by Zakharov \cite{Z}. Exact solutions for the nonlinear evolution
of the unstable line solitons
were obtained and analyzed in \cite{TM93,PS93}. If the transverse perturbation is proportional
to a single unstable eigenmode, it results in the monotonic transformation of the unstable line soliton
with a larger-than-critical speed of propagation
to the travelling solitary wave with the periodic transverse modulation of the same period
and an ejection of a stable
line soliton with a smaller-than-critical speed. In the case of multi-mode perturbations, several modulated
travelling waves and the residual line soliton
are generated in the dynamics of an unstable line soliton, according to
the exact solutions of the KP-I equation \cite{PS93}.

Transformation of the unstable line solitons with a nearly critical speed of propagation was studied
in the framework of the KP-I equation with an asymptotic multi-scale expansion method \cite{GP}. An integrable Eckhaus equation
was derived from the integrable KP-I equation. This Eckhaus equation
correctly represents the monotonic transition of the unstable line soliton
to the transversely modulated solitary wave and a ``radiation" of a stable line soliton
of a smaller speed of propagation, in comparison with the exact solutions to the KP-I equation
(see Section 4.7 in \cite{kivpel}). Similar asymptotic reductions were
reported in \cite{PKS} for the line dark solitons of the defocusing elliptic NLS equation.

The present work is devoted to the justification of the asymptotic model describing the nonlinear
dynamics of the transverse perturbations to the line soliton with a nearly critical speed of propagation.
Stability of these line solitons and transversely modulated solitary waves
are derived as a by-product of this asymptotic model. Unfortunately, the analytical setup does not
apply to the KP-I equation, partly, because the
continuous spectrum of the linearized operators does not move to the left-half plane
in exponentially weighted spaces. For a better model, where this difficulty does not arise,
we consider the Zakharov--Kuznetsov (ZK) equation,
\begin{equation}
\label{ZK}
u_t + 12 u u_x + u_{xxx} + u_{xyy} = 0,
\end{equation}
which features a similar phenomenon. The ZK equation is an anisotropic generalization
of the Korteweg--De Vries (KdV) equation in
two spatial dimensions \cite{ZK}. Justification of the ZK equation in the context of
the Euler--Poisson equations for magnetized plasmas was recently reported in \cite{Lannes}.
Asymptotic stability of two-dimensional solitary waves was considered
in the $L^2$-subcritical ZK equations \cite{Munoz}, which includes the ZK equation (\ref{ZK}).

Transverse stability of line solitons is very similar
between the ZK and KP-I equations \cite{RT0,RT1}, but the ZK equation can be analyzed successfully by
using exponentially weighted Sobolev spaces, similar to the analysis of the KdV equation
\cite{pw1,pw2}. Instability of line solitons in the ZK equation is known for quite some time
in physics literature, e.g., see \cite{Bettison}.

The main result of this work is to derive and to justify the first-order differential equation
\begin{equation}
\label{normal-form-ZK}
\frac{d b}{dt} = \lambda'(c_*) (c_+ - c_*) b + \gamma |b|^2 b, \quad t > 0,
\end{equation}
where $\lambda'(c_*) > 0$, $\gamma < 0$ are real-valued numerical coefficients,
$c_*$ is the critical speed of the line soliton, $c_+ \in \mathbb{R}$ depends on the initial
conditions, and $b(t) : \mathbb{R}_+ \to \mathbb{C}$ is an amplitude
of transverse perturbation. The differential equation (\ref{normal-form-ZK})
describes the nonlinear dynamics of a small transverse perturbation of a fixed period to the line soliton
with a nearly critical speed $c_*$ and is referred to as ``{\em normal form for
transverse instability of the line soliton with a nearly critical speed of propagation}".

Bifurcations and stability of the line solitons under the transverse perturbations of a critical period
were addressed recently by Yamazaki for the elliptic version of the NLS equation in \cite{Y}
and for the ZK equation (\ref{ZK}) in \cite{YY}.

In \cite{Y}, the bifurcation problem is analyzed with the Lyapunov--Schmidt reduction method and
the nonlinear orbital stability is deduced from the energy method based on convexity of the action
functional for the NLS equation, which is the same as for the ZK equation (\ref{ZK}). It is shown in
\cite{Y} that the transversely modulated solitary waves are orbitally stable for the case of
quadratic nonlinearities. Dynamics near such waves was not studied in \cite{Y}.

In \cite{YY}, asymptotic stability of the line solitons with the smaller-than-critical speeds
and the transversely modulated solitary waves with the larger-than-critical speeds
was shown for the ZK equation (\ref{ZK}) with a Liouville-type theorem and virial type estimates.

As the main application of the normal form (\ref{normal-form-ZK}), we show that the conclusions
of \cite{Y,YY} are recovered here with a different technique. In addition, nonlinear dynamics
of perturbations near the line soliton with a nearly critical speed of propagation is clarified
from solutions of the normal form (\ref{normal-form-ZK}).

Although $b(t)$ is a complex amplitude, the normal form (\ref{normal-form-ZK}) describes
motion with the preserved $\arg(b)$, hence dynamics is equivalent to the phase line for $|b|$.
If $c_+ < c_*$, the normal form (\ref{normal-form-ZK}) describes a monotonic exponential decay
$b(t) \to 0$ as $t \to +\infty$ and suggests that the line solitons with the smaller-than-critical
speeds are asymptotically stable with respect to small transverse perturbations of a fixed period,
in agreement with the conclusion of \cite{YY}.
If $c_+ \geq c_*$, the normal form (\ref{normal-form-ZK}) describes a monotonic algebraic
decay $b(t) \to b_*$ as $t \to \infty$, where
\begin{equation}
\label{modulated-wave-ZK}
|b_*| = \sqrt{\frac{\lambda'(c_*) (c_+ - c_*)}{|\gamma|}}, \quad c_+ \geq c_*.
\end{equation}
This solution suggests a transition from the unstable line solitons with the larger-than-critical
speeds to stable transversely modulated solitary waves, which are also asymptotically stable with respect to small
transverse perturbations of the same fixed period, in agreement with the conclusion of \cite{YY}.

The remainder of this paper is organized as follows. Section 2 contains results on transverse instability of line
solitons and bifurcations of transversely modulated solitary waves in the ZK equation (\ref{ZK}).
Formal derivation and justification of the normal form (\ref{normal-form-ZK}) is given in Section 3.\\

{\bf Acknowledgement.} This paper was mostly written in 2014-15 after discussions with F. Rousset
(Paris--Sud) and N. Tzvetkov (Cergy--Pontoise) before papers \cite{Y} and \cite{YY}
were first posted on arXiv and then published. The author thanks the collaborators for useful discussions
and valuable comments.

\section{Transverse instability of line solitons for the ZK equation (\ref{ZK})}

We are concerned here with the transverse instability of the line solitons under periodic
transverse perturbations in the ZK equation (\ref{ZK}). First, we review relevant properties
of the line solitons of the KdV equation. Next, we obtain spectral transverse stability results for
the line solitons with a nearly critical speed of propagation. Further, we study bifurcations of
the transversely modulated solitary waves. Finally, we present the main result on the normal form
for transverse instability of the line solitons with a nearly critical speed of propagation.

\subsection{Properties of line solitons of the KdV equation}

Line solitons are expressed analytically as the two-parameter family
\begin{equation}
\label{line-soliton}
u_c(\xi) = c \;{\rm sech}^2(\sqrt{c} \xi), \quad \xi = x - 4 c t - x_0,
\end{equation}
where $c > 0$ is the speed parameter and $x_0 \in \mathbb{R}$ is the translation parameter.
The line soliton (\ref{line-soliton}) for a fixed $c > 0$ is a critical point of the action functional
\begin{equation}
\label{action}
\Lambda_c(u) := \frac{1}{2} \int_{\mathbb{R}} \left[ \left( \partial_{\xi} u \right)^2 - 4 u^3 + 4c u^2 \right] d \xi,
\end{equation}
which is a linear combination of the energy and momentum of the KdV equation.
The second variation of the action functional (\ref{action}) at the line soliton (\ref{line-soliton})
is defined by the Hessian operator $L_c : H^2(\mathbb{R}) \to L^2(\mathbb{R})$, the differential expression
of which is given by
\begin{equation}
\label{Schrodinger-operator}
L_c := -\partial^2_{\xi} + 4 c - 12 c \; {\rm sech}^2(\sqrt{c} \xi).
\end{equation}
The Schr\"{o}dinger operator $L_c : H^2(\mathbb{R}) \to L^2(\mathbb{R})$ is known \cite{LL} to have the essential
spectrum located on $[4c,\infty)$ and three simple
isolated eigenvalues $\lambda_1 < \lambda_2 < \lambda_3 < 4c$. More precisely,
the eigenvalues and the corresponding non-normalized eigenfunctions of $L_c$ are given by
\begin{eqnarray}
\label{first-eig}
& \lambda_1 = -5c, & \quad \varphi_1 = {\rm sech}^3(\sqrt{c} \xi), \\
\label{second-eig}
& \lambda_2 = 0, & \quad \varphi_2 = {\rm sech}^2(\sqrt{c} \xi) \tanh(\sqrt{c} \xi), \\
& \lambda_3 = 3c, & \quad \varphi_3 = 4 {\rm sech}(\sqrt{c} \xi) - 5 {\rm sech}^3(\sqrt{c} \xi).
\end{eqnarray}
The first two eigenvalues and eigenfunctions are particularly important
in the spectral transverse stability analysis of the line solitons (\ref{line-soliton}).

\subsection{Spectral transverse stability analysis of line solitons}

Let us consider the $2\pi$-periodic transverse perturbation to the line solitons (\ref{line-soliton}).
Therefore, we fix the length of the transverse periodic perturbation and vary the speed parameter
$c$. Substituting $u(x,t) = u_c(\xi) + U(\xi) e^{\lambda t + i k y}$
with $k \in \mathbb{Z}$ into the ZK equation (\ref{ZK}) and dropping the quadratic terms in $U$ yields
the spectral problem
\begin{equation}
\label{spectral}
\partial_{\xi} (L_c + k^2) U = \lambda U, \quad k \in \mathbb{Z},
\end{equation}
where $k$ is the wave number of the $2\pi$-periodic transverse perturbation
and $\partial_{\xi} L_c : H^3(\mathbb{R}) \to L^2(\mathbb{R})$ is the linearized operator for the KdV equation.

For $k = 0$, the spectral problem (\ref{spectral}) coincides with the one for the KdV equation.
It is known from the work of Pego \& Weinstein \cite{pw1} that the spectrum of
$\partial_{\xi} L_c : H^3(\mathbb{R}) \to L^2(\mathbb{R})$
consists of a double zero eigenvalue
and a continuous spectrum on $i \mathbb{R}$. The double zero eigenvalue is associated with the
following Jordan block of the operator $\partial_{\xi} L$:
\begin{equation}
\label{double-zero}
\partial_{\xi} L_c \partial_{\xi} u_c = 0, \quad \partial_{\xi} L_c \partial_c u_c = -4 \partial_{\xi} u_c,
\end{equation}
where the derivatives of $u_c$ in $\xi$ and $c$ are exponentially decaying functions of $\xi$.
The following lemma characterizes the spectral problem (\ref{spectral}) for any $k \in \mathbb{N}$.

\begin{lem}
For any $k \in \mathbb{N}$, the spectral problem (\ref{spectral})
has a pair of real eigenvalues $\pm \lambda_k(c)$
if $c > c_k : = \frac{k^2}{5}$. No eigenvalues with ${\rm Re}(\lambda) \neq 0$
exist if $c \in (0,c_k)$.
\label{theorem-spectral}
\end{lem}

\begin{Proof}
For any $k \in \mathbb{N}$, the self-adjoint Schr\"{o}dinger operator $L_c + k^2 : H^2(\mathbb{R}) \to L^2(\mathbb{R})$
is strictly positive for $c \in (0,c_k)$ and admits a simple negative eigenvalue for $c > c_k$, where $c_k := \frac{k^2}{5}$.
It follows from the main theorem in \cite{KS,Pel} that the spectrum of
$\partial_{\xi} (L_c + k^2) : H^3(\mathbb{R}) \to L^2(\mathbb{R})$ has
exactly one pair of real eigenvalues $\pm \lambda_k(c)$ if $c > c_k$ and no eigenvalues with ${\rm Re}(\lambda) \neq 0$ if $c \in (0,c_k)$.
\end{Proof}

\begin{rem}
By using the energy method in \cite{RT2}, one can actually prove nonlinear orbital stability of
the line solitons with $c \in (0,c_*)$ and nonlinear instability of the line solitons with $c > c_*$, where
$$
c_* := \min\limits_{k \in \mathbb{N}} c_k \equiv \frac{1}{5}.
$$
\end{rem}

The following result gives a precise characterization of the unstable eigenvalue bifurcating at $c = c_*$
in the spectral problem (\ref{spectral}) with $k = 1$. To do so, we introduce
the exponentially weighted space
\begin{equation}
\label{exp-space}
H^s_{\mu}(\mathbb{R}) =
\left\{ u \in H^s_{\rm loc}(\mathbb{R}) : \quad e^{\mu \xi} u \in H^s(\mathbb{R}) \right\}, \quad s \geq 0, \quad \mu > 0.
\end{equation}
This weighted space is used to push the continuous spectrum of the operator
$\partial_{\xi} L_c : H^3_{\mu}(\mathbb{R}) \to L^2_{\mu}(\mathbb{R})$
to the left-half plane for $\mu > 0$ sufficiently small \cite{pw1,pw2}.
In the generalized KdV equation with a triple zero eigenvalue, the exponentially weighted space is
used to construct the Jordan block for the triple zero eigenvalue in \cite{CCP,Marzuola}.

\begin{lem}
There is $\mu_0 > 0$ such that for every $\mu \in (0,\mu_0)$, the spectral problem (\ref{spectral})
with $k = 1$ and $c = c_*$ considered in $L^2_{\mu}(\mathbb{R})$ admits a simple zero eigenvalue with the eigenfunction
$\psi_* \in H^3_{\mu}(\mathbb{R})$ and the adjoint eigenfunction $\eta_* \in H^3_{-\mu}(\mathbb{R})$, where
\begin{equation}
\label{limiting-function}
\psi_*(\xi) = {\rm sech}^3(\sqrt{c_*} \xi), \quad
\eta_*(\xi) = \int_{-\infty}^{\xi} {\rm sech}^3(\sqrt{c_*} \xi') d\xi',
\end{equation}
hence,
\begin{equation}
\label{limiting-coefficient}
\langle \eta_*, \psi_* \rangle_{L^2} = \frac{1}{2} \left( \int_{\mathbb{R}} {\rm sech}^3(\sqrt{c_*} \xi) d \xi \right)^2
= \frac{\pi^2}{8 c_*}.
\end{equation}
Moreover, for a given $\mu \in (0,\mu_0)$, there exist an interval $(c_-,c_+)$ with
 $c_- < c_* < c_+$ such that the spectral problem (\ref{spectral})
with $k = 1$ and $c \in (c_-,c_+)$ considered in $L^2_{\mu}(\mathbb{R})$
admits a small eigenvalue $\lambda(c)$, where the mapping $c \mapsto \lambda$ is smooth and
is given by
\begin{equation}
\label{lambda-expansion}
\lambda(c) = \lambda'(c_*) (c-c_*) + \mathcal{O}((c-c_*)^2) \quad \mbox{\rm as} \quad c \to c_*,
\end{equation}
with
\begin{equation}
\label{derivative-eigenvalue}
\lambda'(c_*) = \frac{128}{3 \pi^2} \sqrt{c_*}.
\end{equation}
\label{theorem-bifurcation}
\end{lem}

\begin{Proof}
Existence of the zero eigenvalue of $\partial_{\xi} (L_{c_*} + 1) : H^3(\mathbb{R}) \to L^2(\mathbb{R})$
follows from existence of the negative eigenvalue $-1$ of $L_{c_*}$ in (\ref{first-eig}).
The eigenfunction $\psi_* = \varphi_1$ at $c = c_*$ is exponentially decaying in $\xi$,
therefore, there is a positive $\mu_0$ such that $e^{\mu \xi} \psi_*$ is exponentially
decaying in $\xi$ for every $\mu \in (0,\mu_0)$.  Thus, the operator
$\partial_{\xi} (L_{c_*}+1) : H^3_{\mu}(\mathbb{R}) \to L^2_{\mu}(\mathbb{R})$  has a simple zero eigenvalue.

Let us construct an adjoint operator to the operator
$\partial_{\xi} (L_{c_*}+1) : H^3_{\mu}(\mathbb{R}) \to L^2_{\mu}(\mathbb{R})$.
We take a standard $L^2$ product with the differential expression
$e^{\mu \xi} \partial_{\xi} (L_{c_*}+1) e^{-\mu \xi}$, integrate by parts,
and obtain the adjoint differential expressions
in the form $-e^{-\mu \xi} (L_{c_*}+1) \partial_{\xi} e^{\mu \xi}$.
Replacing exponential factors by weighted spaces yields the adjoint operator in
the form $-(L_{c_*}+1)\partial_{\xi} : H^3_{-\mu}(\mathbb{R}) \to L^2_{-\mu}(\mathbb{R})$.

For $\mu > 0$, the eigenfunction $\eta_*$ of the kernel of
$-(L_{c_*}+1) \partial_{\xi} : H^3_{-\mu}(\mathbb{R}) \to L^2_{-\mu}(\mathbb{R})$
must decay faster than $e^{-\mu \xi}$ grows as $\xi \to -\infty$.
Since $\partial_{\xi} \eta_*$ is proportional to $\psi_*$ and $\eta_*(\xi) \to 0$ as $\xi \to -\infty$,
we set
$$
\eta_*(\xi) := \int_{-\infty}^{\xi} \psi_*(\xi') d \xi'
$$
and obtain (\ref{limiting-function}).
The nonzero inner product in (\ref{limiting-coefficient}) ensures simplicity of the Jordan block
for the zero eigenvalue of $\partial_{\xi} (L_{c_*}+1) : H^3_{\mu}(\mathbb{R}) \to L^2_{\mu}(\mathbb{R})$
with $\mu \in (0,\mu_0)$.

It remains to prove the last assertion of the lemma with the asymptotic expansion (\ref{lambda-expansion}).
Since the simple zero eigenvalue of $\partial_{\xi} (L_c + 1) : H^3_{\mu}(\mathbb{R}) \to L^2_{\mu}(\mathbb{R})$ at $c = c_*$
is isolated from the continuous spectrum of this operator for any fixed $\mu \in (0,\mu_0)$,
one can use analytic perturbation theory \cite{Kato}. In particular,
the operator $L_c : H^2(\mathbb{R}) \mapsto L^2(\mathbb{R})$ is analytic at $c = c_*$
and admits the following Taylor expansion
\begin{equation}
\label{expansion-operator}
L_c = L_{c_*} + L_{c_*}' (c-c_*) + \tilde{L}_c (c - c_*)^2,
\end{equation}
where
\begin{equation}
\label{derivative-operator}
L_{c_*}' = 4 - 12 {\rm sech}^2(\sqrt{c_*}\xi) + 12 \sqrt{c_*} \xi \; {\rm sech}^2(\sqrt{c_*} \xi) \; \tanh(\sqrt{c_*} \xi)
\end{equation}
and $\tilde{L}_{c}$ is an exponentially decaying and bounded potential as $c \to c_*$.
By using formal expansions,
\begin{equation}
\label{expansion-lambda}
\lambda(c) = \lambda_1 (c-c_*) + \mathcal{O}((c-c_*)^2), \quad \psi(c) = \psi_* + \psi_1 (c-c_*)
+ \mathcal{O}_{H^3_{\mu}}((c-c_*)^2),
\end{equation}
we obtain the linear inhomogeneous equation at the order of $\mathcal{O}(c-c_*)$:
\begin{equation}
\label{first-order-perturbation-theory}
\partial_{\xi} (L_{c_*} + 1) \psi_1 + \partial_{\xi} L_{c_*}' \psi_* = \lambda_1 \psi_*.
\end{equation}
This equation is considered in $L^2_{\mu}(\mathbb{R})$ for a fixed $\mu \in (0,\mu_0)$.
Projecting (\ref{first-order-perturbation-theory}) to $\eta_*$,
the eigenfunction in the kernel of the adjoint operator
$-(L_{c_*}+1)\partial_{\xi} : H^3_{-\mu}(\mathbb{R}) \mapsto L^2_{-\mu}(\mathbb{R})$, we obtain
\begin{equation}
\label{lambda-1-explicit}
\lambda_1 \langle \eta_*,\psi_* \rangle_{L^2} = \langle \eta_*, \partial_{\xi} L_{c_*}' \psi_* \rangle_{L^2} =
- \langle \psi_*, L_{c_*}' \psi_* \rangle_{L^2}.
\end{equation}
From (\ref{limiting-coefficient}), (\ref{derivative-operator}), and (\ref{lambda-1-explicit}), we obtain
$$
\langle \psi_*, L_{c_*}' \psi_* \rangle_{L^2} = -\frac{16}{3 \sqrt{c_*}} \quad \Rightarrow \quad
\lambda_1 = \frac{128}{3 \pi^2} \sqrt{c_*},
$$
which agrees with (\ref{derivative-eigenvalue}) because $\lambda_1 = \lambda'(c_*)$. Justification
of the asymptotic expansion (\ref{expansion-lambda}) is developed with the analytic perturbation theory
(Theorem 1.7 in \cite[Chapter VII]{Kato}).
\end{Proof}

\begin{rem}
Since the continuous spectrum of $\partial_{\xi}(L_c + 1) : H^3_{\mu}(\mathbb{R}) \to L^2_{\mu}(\mathbb{R})$
is located in the left-half plane for $\mu > 0$ \cite{pw1}, whereas the eigenvalue $\lambda(c)$ is positive for $c > c_*$,
we can send $\mu \to 0^+$ without affecting the eigenvalue $\lambda(c) > 0$. Therefore, the eigenvalue $\lambda(c)$
persists for the operator $\partial_{\xi}(L_c + 1) : H^3(\mathbb{R}) \to L^2(\mathbb{R})$.
By symmetry $\xi \mapsto -\xi$ and $\lambda \mapsto -\lambda$, the operator
$\partial_{\xi}(L_c + 1) : H^3(\mathbb{R}) \to L^2(\mathbb{R})$ also admits the eigenvalue
$-\lambda(c) < 0$ for the same case $c > c_*$.
Thus, the spectral problem (\ref{spectral}) with $k = 1$ has a pair of real eigenvalues $\pm \lambda(c)$
for $c > c_*$, in agreement with Lemma \ref{theorem-spectral}.
\end{rem}

\begin{rem}
Since the eigenvalue $\lambda(c)$ in $L^2_{\mu}(\mathbb{R})$
is negative for $c < c_*$, we cannot send $\mu \to 0^+$ without affecting the
eigenvalue $\lambda(c) < 0$ by the deformation of the continuous spectrum of
$\partial_{\xi}(L_c+1) : H^3_{\mu}(\mathbb{R}) \to L^2_{\mu}(\mathbb{R})$.
As is well-known \cite{pw1}, the eigenvalue $\lambda(c) < 0$
of the operator $\partial_{\xi} (L_c + 1) : H^3_{\mu}(\mathbb{R}) \to L^2_{\mu}(\mathbb{R})$
becomes a resonant pole of the operator $\partial_{\xi} (L_c + 1) : H^3(\mathbb{R}) \to L^2(\mathbb{R})$ for $c < c_*$.
\end{rem}

\begin{rem}
Perturbation expansions (\ref{expansion-operator}) and (\ref{expansion-lambda}) can also be used
to characterize the shift of the zero eigenvalue of the operator $L_c + 1$ in $L^2(\mathbb{R})$ for $c \neq c_*$. In this case,
$\psi_1$ satisfies
$$
(L_{c_*}+1) \psi_1 + L'_{c_*} \psi_* = \lambda_1 \psi_1,
$$
from which we obtain
$$
\lambda_1 \| \psi_* \|_{L^2}^2 = \langle \psi_*, L_{c_*}' \psi_* \rangle_{L^2}.
$$
Since $\langle \psi_*, L_{c_*}' \psi_* \rangle_{L^2} < 0$, we have $\lambda_1 < 0$.
Therefore, the zero eigenvalue of the operator $L_{c_*} + 1$ in $L^2(\mathbb{R})$
becomes a negative eigenvalue of this operator for $c > c_*$. This shift induces the spectral
instability of the line soliton, in agreement with Lemma \ref{theorem-spectral}.
\end{rem}

\subsection{Transversely modulated solitary waves}

The instability bifurcation of the line solitons in Lemma \ref{theorem-bifurcation} is related to the bifurcation
of a new family of travelling solitary waves with the periodic transverse modulation.
Such transversely modulated solitary waves satisfy the nonlinear elliptic problem
\begin{equation}
\label{elliptic-problem}
-\frac{\partial^2 u}{\partial \xi^2} - \frac{\partial^2 u}{\partial y^2} + 4 c u - 6 u^2 = 0, \quad
(x,y) \in \mathbb{R} \times \mathbb{T},
\end{equation}
where $\mathbb{T}$ is a $2\pi$-periodic torus.
In order to eliminate the translational symmetries in $\xi$ and $y$, we define the space of even functions
both in $\xi$ and $y$:
\begin{equation}
\label{even-space}
H^s_{\rm even} = \left\{ u \in H^s(\mathbb{R} \times \mathbb{T}) : \quad u(-\xi,y) = u(\xi,y) = u(\xi,-y) \right\}, \quad s \geq 0.
\end{equation}
The following lemma describes bifurcation of the transversely modulated solitary waves
from the line soliton with a nearly critical speed of propagation.

\begin{lem}
\label{lemma-modulated-wave}
There exists $c_+ > c_*$ such that for every $c \in (c_*,c_+)$,
the nonlinear elliptic problem (\ref{elliptic-problem}) has a nontrivial solution $u_b$ in $H^2_{\rm even}$
in addition to the line soliton (\ref{line-soliton}). The solution $u_b$ is expressed by
the expansion
\begin{equation}
\label{LS-decomposition}
u_b(\xi,y) = u_{c_*}(\xi) + 2 b \cos(y) \psi_*(\xi) + \tilde{u}_b(\xi,y),
\end{equation}
where $b \in \mathbb{R}$ is a nonzero root of the algebraic equation
\begin{equation}
\label{normal-form-static}
\alpha (c-c_*) b + \beta b^3 = 0
\end{equation}
and $\tilde{u}_b \in H^2_{\rm even}$ satisfies the bound $\| \tilde{u}_b \|_{H^2} \leq A b^2$
for a positive constant $A$ independently of $b$ and $c$. Here
$$
\alpha = -\langle \psi_*, L_{c_*}' \psi_* \rangle_{L^2} = \frac{16}{3 \sqrt{c_*}} > 0
$$
and $\beta < 0$ is a numerical coefficient given by (\ref{gamma-stationary}) below.
\end{lem}

\begin{Proof}
The proof is close to the Crandall-Rabinowitz local bifurcation theory \cite{CR71} and
relies on the method of Lyapunov--Schmidt reduction \cite{Nirenberg}.

We write the decomposition (\ref{LS-decomposition}) in $H^2_{\rm even}$
equipped with the orthogonality condition $\langle v_1, \tilde{u}_b \rangle_{L^2(\mathbb{R} \times \mathbb{T})} = 0$,
where $v_1(\xi,y) := \cos(y) \psi_*(\xi)$ is the eigenfunction of the kernel
of $L_{c_*} - \partial_y^2 : H^2_{\rm even} \to L^2_{\rm even}$.
Let $\Pi$ be an orthogonal projection operator in $L^2(\mathbb{R} \times \mathbb{T})$
in the complement of $v_1$. For every $c = c_* + \delta$ with $\delta \in \mathbb{R}$ being
sufficiently small, the correction term $\tilde{u}_b \in H^2_{\rm even}$ and
the parameter $b \in \mathbb{R}$ are defined by the projection equations
of the Lyapunov--Schmidt reduction method:
\begin{equation}
\label{reduction-1}
(L_{c_*} - \partial_y^2 + 4 \delta) \tilde{u}_b = \Pi \tilde{F},
\end{equation}
with
\begin{equation}
\label{expression-F}
\tilde{F} := - 4 \delta u_{c_*} - 8 \delta b \cos(y) \psi_* + 6 (2 b \cos(y) \psi_* + \tilde{u}_b)^2
\end{equation}
and
$$
\frac{1}{2\pi} \langle v_1, \tilde{F} \rangle_{L^2(\mathbb{R} \times \mathbb{T})} = 0,
$$
which yields
\begin{equation}
\label{reduction-2} - 4 \delta b \| \psi_* \|^2_{L^2(\mathbb{R})} +
\frac{3}{\pi} \langle \cos(y) \psi_*, (2 b \cos(y) \psi_* + \tilde{u}_b)^2 \rangle_{L^2(\mathbb{R} \times \mathbb{T})} = 0.
\end{equation}
Thanks to the symmetry constraints in (\ref{even-space}), the kernel of the linear operator
$L_{c_*} - \partial_y^2 : H^2_{\rm even} \to L^2_{\rm even}$ is one-dimensional. By the spectral calculus,
there are $B > 0$ and $\delta_0 > 0$ such that
\begin{equation}
\label{inverse-operator}
\| \Pi (L_{c_*} - \partial_y^2 + 4 \delta)^{-1} \Pi \|_{L^2_{\rm even} \to L^2_{\rm even}} \leq B, \quad \forall |\delta| < \delta_0.
\end{equation}
Thanks to (\ref{inverse-operator}), for every small $b \in \mathbb{R}$ and small $\delta \in \mathbb{R}$,
the fixed-point argument can be applied to
solve equation (\ref{reduction-1}) with (\ref{expression-F}) in $H^2_{\rm even}$
and to obtain a unique $\tilde{u}_b \in H^2_{\rm even}$ satisfying the bound
\begin{equation}
\label{reduction-2a}
\| \tilde{u}_b \|_{H^2} \leq A (|\delta| + b^2),
\end{equation}
where the positive constant $A$ is independent of $\delta$ and $b$. However, the solution $\tilde{u}_b \in H^2_{\rm even}$ with
the bound (\ref{reduction-2a})  is insufficient for the derivation of the algebraic equation (\ref{normal-form-static})
from the bifurcation equation (\ref{reduction-2}) \cite{Phan}.

In order to obtain the algebraic equation (\ref{normal-form-static}), we perform a near-identity transformation
\begin{equation}
\label{v-component}
\tilde{u}_b(\xi,y) = 2 b^2 \cos(2y) w_2(\xi) + b^2 w_0(\xi) + \delta \partial_c u_{c_*}(\xi) + \tilde{w}(\xi,y),
\end{equation}
where
\begin{equation}
\label{explicit-solution-2}
\partial_c u_{c_*}(\xi) = {\rm sech}^2(\sqrt{c_*} \xi) - \sqrt{c_*} \xi \tanh(\sqrt{c_*} \xi) {\rm sech}^2(\sqrt{c_*} \xi),
\end{equation}
the correction terms $w_0$ and $w_2$ are found from the linear inhomogeneous equations:
\begin{equation}
\label{inhom-eq-1}
L_{c_*} w_0 = 12 \psi_*^2
\end{equation}
and
\begin{equation}
\label{inhom-eq-2}
(L_{c_*} + 4) w_2 = 6 \psi_*^2,
\end{equation}
whereas $\tilde{w}$ satisfies the transformed equation
\begin{eqnarray}
\label{reduction-3}
(L_{c_*} - \partial_y^2 + 4 \delta) \tilde{w} = \Pi \tilde{G},
\end{eqnarray}
with
\begin{eqnarray*}
\tilde{G} :=  - 4 \delta^2 \partial_c u_{c_*}
- 8 \delta b^2 \cos(2y) w_2 - 4 \delta b^2 w_0
- 8 \delta b \cos(y) \psi_* + 24 b \cos(y) \psi_* \tilde{u}_b + 6 \tilde{u}_b^2,
\end{eqnarray*}
where $\tilde{u}_b$ is related to $\tilde{w}$ by (\ref{v-component}).
By the same argument as above, for every small $b \in \mathbb{R}$ and $\delta \in \mathbb{R}$,
there exists a unique solution $\tilde{w} \in H^2_{\rm even}$ of equation (\ref{reduction-3})
satisfying the bound
\begin{equation}
\label{reduction-3a}
\| \tilde{w} \|_{H^2} \leq A (\delta^2 + |\delta| |b| + |b|^3),
\end{equation}
where the positive constant $A$ is independent of $\delta$ and $b$.

Substituting the near-identity transformation (\ref{v-component}) into the bifurcation equation
(\ref{reduction-2}) and using the bound (\ref{reduction-3a})
for the component $\tilde{w} \in H^2_{\rm even}$, we rewrite
(\ref{reduction-2}) in the equivalent form
\begin{equation}
\label{reduction-4}
\alpha \delta b + \beta b^3 + \mathcal{O}(\delta b^2,b^4) = 0,
\end{equation}
where we have introduced numerical coefficients $\alpha$ and $\beta$ as follows:
\begin{equation}
\label{alpha-static}
\alpha := -4 \| \psi_* \|^2_{L^2} + 12 \langle \psi_*^2, \partial_c u_{c_*} \rangle_{L^2}
\end{equation}
and
\begin{equation}
\label{gamma-static}
\beta := 12 \langle \psi_*^2, w_0 + w_2 \rangle_{L^2}.
\end{equation}
We have also removed $\mathcal{O}(\delta^2)$ from the remainder term in (\ref{reduction-4}), because
the bifurcation equation (\ref{reduction-2}) is identically satisfied
in the case of the line soliton with $b = 0$ for every small $\delta \in \mathbb{R}$.

Comparison of (\ref{explicit-solution-2}) and
(\ref{alpha-static}) with (\ref{derivative-operator}) shows that
$$
\alpha = -\langle \psi_*, L_{c_*}' \psi_* \rangle_{L^2} = \frac{16}{3 \sqrt{c}_*}.
$$
On the other hand, the coefficient $\beta$ in (\ref{gamma-static}) is less explicit.
In order to show that $\beta < 0$, we obtain the unique even solution of
the linear inhomogeneous equation (\ref{inhom-eq-1}) in the explicit form
\begin{equation}
\label{explicit-solution-1}
w_0(\xi) = -15 {\rm sech}^2(\sqrt{c_*} \xi) + \frac{15}{2} {\rm sech}^4(\sqrt{c_*}\xi),
\end{equation}
where we have used $c_* = \frac{1}{5}$.
Because $L_{c_*} + 4 : H^2(\mathbb{R}) \to L^2(\mathbb{R})$ is strictly positive, there exists a unique solution
of the linear inhomogeneous equation (\ref{inhom-eq-2}).
Unfortunately, it is not available in the explicit form. Nevertheless, we can represent
this unique solution of equation (\ref{inhom-eq-2}) in the form
\begin{equation}
\label{explicit-solution-3}
w_2(\xi) = 5 {\rm sech}^2(\sqrt{c_*}\xi) + \frac{15}{4} {\rm sech}^4(\sqrt{c_*} \xi) - \tilde{w}_2(\xi),
\end{equation}
where $\tilde{w}_2$ is found from the inhomogeneous equation
\begin{equation}
\label{inhomogeneous-eq-tilde-w-2}
(L_{c_*}+4) \tilde{w}_2(\xi) = 20 {\rm sech}^2(\sqrt{c_*} \xi).
\end{equation}
By the maximum principle for the elliptic operator $(L_{c_*} + 4)$,
the component $\tilde{w}_2$ satisfies $\tilde{w}_2(\xi) \geq 0$ for all $\xi \in \mathbb{R}$.
After computing the integrals in (\ref{gamma-static}), we obtain
\begin{equation}
\label{gamma-stationary}
\beta = 12 \langle \psi_*^2, w_0 + w_2 \rangle_{L^2} = -12 \langle \psi_*^2, \tilde{w}_2 \rangle_{L^2} < 0.
\end{equation}
The cubic equation (\ref{normal-form-static}) follows from the truncation of the bifurcation
equation (\ref{reduction-4}) at the first two terms. This is justified because if $b \neq 0$,
then it follows from (\ref{reduction-4}) that
$$
\delta = -\frac{\beta}{\alpha} b^2 + \mathcal{O}(b^3)
$$
and if $b_0$ is a root of the cubic equation (\ref{normal-form-static}) for a given small
$\delta = c - c_*$, then
$$
|b-b_0| \leq A b_0^2,
$$
where the positive constant $A$ is independent of $\delta$. Then, it follows from the bound
(\ref{reduction-2a}) that $\| \tilde{u}_b \|_{H^2} \leq A b_0^2$, which justifies
the decomposition (\ref{LS-decomposition}) after the change of notation $b_0 \mapsto b$.
\end{Proof}

\begin{rem}
Since $\alpha > 0$ and $\beta < 0$, the nonzero solutions for $b$ exists in the cubic equation
(\ref{normal-form-static}) if and only if $c > c_*$. By Lemma \ref{theorem-spectral},
the line soliton with $c > c_*$
is spectrally unstable with respect to the transverse $2\pi$-periodic perturbations.
\end{rem}

\subsection{Statement of the main theorem}

The ZK equation (\ref{ZK}) was shown in \cite{Linares} to be
locally well-posed in $H^s(\mathbb{R} \times \mathbb{T})$ for $s > \frac{3}{2}$ and globally well-posed
for perturbations of the line solitons (\ref{line-soliton}) in $H^1(\mathbb{R}^2)$. More recently,
the ZK equation (\ref{ZK}) was shown in  \cite{Moliner} to be globally well-posed in $H^1(\mathbb{R} \times \mathbb{T})$.
The latter well-posedness result allows us to employ the energy method in the justification of the normal form
for transverse instability of the line soliton with a nearly critical speed of propagation.

Let us denote by $H^s_{\mu}(\mathbb{R} \times \mathbb{T})$ the exponentially weighted version
of the space $H^s(\mathbb{R} \times \mathbb{T})$, $s \geq 0$ with the weight $\mu > 0$ applied
in the $\xi$ axis only.

In order to characterize the dynamics of solutions of the ZK equation (\ref{ZK}) in time $t$
near the line solitons (\ref{line-soliton}), we introduce varying parameters $a(t)$ and
$c(t)$ of the line solitons as well as its perturbation $\tilde{u}(t)$ defined
in $H^1(\mathbb{R} \times \mathbb{T}) \cap H^1_{\mu}(\mathbb{R} \times \mathbb{T})$ for every $t \in \mathbb{R}_+$.
Hence, we introduce the travelling coordinate $\xi = x - 4 a(t)$ and use the decomposition
\begin{equation}
\label{decomposition}
u(x,y,t) = u_{c(t)}(\xi) + \tilde{u}(\xi,y,t), \quad \xi = x - 4 a(t).
\end{equation}
The time evolution of the varying parameters $a(t)$ and $c(t)$ and the perturbation term $\tilde{u}(t)$
are to be found from the evolution problem
\begin{equation}
\label{amplitude-evolution}
\tilde{u}_t = \partial_{\xi} (L_{c} - \partial_y^2 + 4(\dot{a}-c)) \tilde{u}
+ 4 (\dot{a} - c) \partial_{\xi} u_{c} - \dot{c} \partial_c u_{c}
- 6 \partial_{\xi} \tilde{u}^2,
\end{equation}
where the differential expression for $L_c$ is given by (\ref{Schrodinger-operator}).

Although parameters $a(t)$ and $c(t)$ vary in time along the time evolution of system (\ref{amplitude-evolution}),
we will prove that $\dot{a}(t)$ remains close to $c(t)$ and that $c(t)$ remains close to $c_*$ for all times,
where $c_*$ is the critical speed
of propagation, see bounds (\ref{final-bound}) below. Therefore, the representation (\ref{decomposition})
is used to characterize dynamics of the line solitons with a nearly critical speed of propagation.

For $c = c_*$, the operators $L_{c_*} + k^2 : H^2(\mathbb{R}) \to L^2(\mathbb{R})$ are coercive
for any $k \in \mathbb{Z} \backslash \{0,\pm 1\}$.
This property will be used to control perturbations with the corresponding Fourier wave numbers.
On the other hand, the operator $\partial_{\xi} L_{c_*} : H^3_{\mu}(\mathbb{R}) \to L^2_{\mu}(\mathbb{R})$ has a double
zero eigenvalue associated with the Jordan block (\ref{double-zero}), whereas
the operator $\partial_{\xi} (L_{c_*} + 1)  : H^3_{\mu}(\mathbb{R}) \to L^2_{\mu}(\mathbb{R})$ has a simple zero eigenvalue by
Lemma \ref{theorem-bifurcation}. The double zero eigenvalue of $\partial_{\xi} L_{c_*}  : H^3_{\mu}(\mathbb{R}) \to L^2_{\mu}(\mathbb{R})$
is already incorporated in the decomposition (\ref{decomposition}), whereas
the simple zero eigenvalue of $\partial_{\xi} (L_{c_*}+1)  : H^3_{\mu}(\mathbb{R}) \to L^2_{\mu}(\mathbb{R})$ will be
incorporated in the secondary decomposition of the perturbation term $\tilde{u}$.

The following theorem represents the normal form for transverse instability of
the line soliton with a nearly critical speed of propagation.

\begin{theo}
Consider the Cauchy problem for the evolution equation (\ref{amplitude-evolution}) with
$$
\tilde{u}(0) \in H^1(\mathbb{R} \times \mathbb{T})\cap H^1_{\mu}(\mathbb{R} \times \mathbb{T}),
$$
where $\mu > 0$ is sufficiently small. There exist $\eps_0 > 0$ and $C_0 > 0$ such that if the initial data satisfy the bound
\begin{equation}
\label{initial-bound}
\| \tilde{u}(0) - 2 \eps \cos(y) \psi_* \|_{H^1(\mathbb{R} \times \mathbb{T}) \cap H^1_{\mu}(\mathbb{R} \times \mathbb{T}) }
+ |c(0) - c_*| \leq \eps^2,
\end{equation}
for every $\eps \in (0,\eps_0)$, then there exist unique functions $a, b, c \in C^1(\mathbb{R}_+)$
and the unique solution
$$
\tilde{u}(t) \in C(\mathbb{R}_+;H^1(\mathbb{R} \times \mathbb{T})\cap H^1_{\mu}(\mathbb{R} \times \mathbb{T}))
$$
of the evolution equation (\ref{amplitude-evolution}) satisfying the bound
\begin{equation}
\label{final-bound}
\| \tilde{u}(t) - (b(t) e^{iy} + \bar{b}(t) e^{-iy}) \psi_* \|_{H^1(\mathbb{R} \times \mathbb{T})}
+ |c(t) - c_*| + |\dot{a}(t) - c(t) | \leq C_0 \eps^2, \quad t \in \mathbb{R}_+,
\end{equation}
Furthermore, the function $b(t)$ satisfies the normal form
\begin{equation}
\label{normal-form}
\dot{b} = \lambda'(c_*) (c_+ - c_*) b + \gamma |b|^2 b, \quad t \in \mathbb{R}_+,
\end{equation}
with $b(0) = \eps$ and $c_+ \in \mathbb{R}$ satisfying $|c_+-c_*| \leq C_0 \eps^2$,
where $\lambda'(c_*) > 0$ is given by (\ref{derivative-eigenvalue}) and
$\gamma < 0$ is a specific numerical coefficient given by (\ref{gamma}) below.
Consequently, $|b(t)| \leq C_0 \eps$ for every $t \in \mathbb{R}_+$.
\label{theorem-normal-form}
\end{theo}

\begin{rem}
The linear part of the normal form (\ref{normal-form}) reproduces the
spectral stability result of Lemma \ref{theorem-bifurcation}.
\end{rem}

\begin{rem}
The stationary part of the normal form (\ref{normal-form})
represents the bifurcation result (\ref{normal-form-static}) in Lemma \ref{lemma-modulated-wave},
except for the quantitative discrepancy between the numerical coefficients
$(\alpha,\beta)$ and numerical coefficients $(\lambda'(c_*),\gamma)$,
see Remark \ref{remark-coefficients} below.
\end{rem}

\begin{rem}
The important fact $\gamma < 0$ is only established with a numerical computation,
see Remark \ref{remark-numerics} below.
\end{rem}

\begin{rem}
We will obtain in the expansions (\ref{decomposition-a-c}) and (\ref{normal-form-1e}) below that
the constant $c_+ - c_*$ is related to the initial data $c(0)$ and $b(0) = \eps$ as follows:
$$
c_+ - c_* = c(0) - c_* + \frac{16}{3} \eps^2 + \mathcal{O}(\eps^4).
$$
Hence, the sign of $c_+ - c_*$, which is crucial to distinguish two different solutions of the normal form
(\ref{normal-form}), depends on the initial data. In particular, if $c(0) > c_*$, then $c_+ > c_*$
and the line soliton is unstable with the perturbation growing towards the stable transversely modulated
solitary waves.
\end{rem}

\begin{rem}
Similarly to the asymptotic stability result in \cite{pw2}, we anticipate that for $\mu > 0$ sufficiently small,
there exists $b_{\infty}$ and $c_{\infty}$ such that the solution in Theorem \ref{theorem-normal-form} satisfies the
following limits:
\begin{equation}
\label{final-bound-asymptotic-smaller}
c_+ \leq c_* : \quad  \lim_{t \to \infty} \| u(t) - u_{c_{\infty}} \|_{H^1_{\mu}(\mathbb{R} \times \mathbb{T})} = 0,
\end{equation}
and
\begin{equation}
\label{final-bound-asymptotic-larger}
c_+ > c_* : \quad \lim_{t \to \infty} \| u(t) - u_{b_{\infty}} \|_{H^1_{\mu}(\mathbb{R} \times \mathbb{T})} = 0,
\end{equation}
where $u_{c_{\infty}}$ is the line soliton (\ref{line-soliton}) with $c = c_{\infty}$ and
$u_{b_{\infty}}$ is the transversely modulated solitary wave defined by Lemma \ref{theorem-bifurcation}
with $b = b_{\infty}$.
The limits (\ref{final-bound-asymptotic-smaller}) and (\ref{final-bound-asymptotic-larger}) are in agreement with
the asymptotic stability results obtained in \cite{YY}. However, the proof of these limits requires more
control of the modulation equations for perturbations in $H^1_{\mu}(\mathbb{R} \times \mathbb{T})$.
The relevant tools are not available from the previous work \cite{pw2}, where
perturbations were considered in $H^1_{\mu}(\mathbb{R})$.
\end{rem}

\section{Derivation and justification of the normal form (\ref{normal-form})}

We first derive the general modulation equations for the varying parameters $a$ and $c$
in the decomposition (\ref{decomposition}) near the line solitons (\ref{line-soliton}).
Next, we introduce the varying parameter $b$ in the secondary decomposition along the neutral eigenmode
$\psi_*$ in Lemma \ref{theorem-bifurcation} and derive the corresponding modulation equation for $b$.
Further, we justify the bound (\ref{final-bound}) with the energy method.
Finally, we simplify the modulation equations and derive the normal form (\ref{normal-form})
by means of nearly identity transformations and the momentum conservation.

\subsection{Modulation equations for parameters $a$ and $c$}

Let $X_c = {\rm span}\{\partial_{\xi} u_c, \partial_c u_c \}$ be an invariant subspace of $L^2(\mathbb{R})$ for the double zero
eigenvalue of the linearized operator $\partial_{\xi} L_c : H^3(\mathbb{R}) \to L^2(\mathbb{R})$,
according to the Jordan block (\ref{double-zero}). Thanks to the exponential decay of $u_c(\xi)$ as $|\xi| \to \infty$,
there is $\mu_0 > 0$ such that $X_c$ is also an invariant subspace of $L^2_{\mu}(\mathbb{R})$ for
$\partial_{\xi} L_c : H^3_{\mu}(\mathbb{R}) \to L^2_{\mu}(\mathbb{R})$ for $\mu \in (0,\mu_0)$.
Similarly, $X_c^* = {\rm span}\{ u_c, \partial_{\xi}^{-1} \partial_c u_c \}$ is an invariant subspace
of $L^2_{-\mu}(\mathbb{R})$ for the double zero eigenvalue of the adjoint operator
$-L_c \partial_{\xi} : H^3_{-\mu}(\mathbb{R}) \to L^2_{-\mu}(\mathbb{R})$.
Recall that the exponentially weighted space is defined in (\ref{exp-space})
and
$$
\partial_{\xi}^{-1} u(\xi) := \int_{-\infty}^{\xi} u(\xi') d \xi'.
$$
In order to avoid confusion between spaces $L^2_{\mu}(\mathbb{R})$ and $L^2_{\mu}(\mathbb{R} \times \mathbb{T})$,
we specify whether the spatial domain is $\mathbb{R}$ or $\mathbb{R} \times \mathbb{T}$ in the $L^2$ inner products and
their induced norms.

The following lemma states the validity of the decomposition (\ref{decomposition}).

\begin{lem}
\label{lem-decomposition}
There exists $\eps_0 > 0$, $\mu_0 > 0$, and $C_0 > 0$ such that if
$u \in C(\mathbb{R}_+,H^1(\mathbb{R} \times \mathbb{T}) \cap H^1_{\mu}(\mathbb{R} \times \mathbb{T}))$
with $\mu \in (0,\mu_0)$ is a global solution to the ZK equation (\ref{ZK}) satisfying
\begin{equation}
\label{orbit-bound}
\eps := \inf_{a \in \mathbb{R}} \| u(x+4a,y,t) - u_{c_*}(x) \|_{H^1(\mathbb{R} \times \mathbb{T}) \cap H^1_{\mu}(\mathbb{R} \times \mathbb{T})}
\leq \eps_0, \quad t \in \mathbb{R}_+,
\end{equation}
then there exist $a, c \in C(\mathbb{R}_+)$ and $\tilde{u} \in
C(\mathbb{R}_+,H^1(\mathbb{R} \times \mathbb{T}) \cap H^1_{\mu}(\mathbb{R} \times \mathbb{T}))$
such that the decomposition
\begin{equation}
\label{decomposition-lemma}
u(x,y,t) = u_{c(t)}(\xi) + \tilde{u}(\xi,y,t), \quad \xi = x - 4 a(t)
\end{equation}
holds with $\tilde{u}(t) \in [X_{c(t)}^*]^{\perp}$ for every $t \in \mathbb{R}_+$, where
\begin{equation}
\label{symplectic-orthogonality}
[X_{c(t)}^*]^{\perp} = \left\{ \tilde{u} \in L^2_{\mu}(\mathbb{R} \times \mathbb{T}) : \quad
\langle u_{c(t)}, \tilde{u} \rangle_{L^2(\mathbb{R} \times \mathbb{T})} =
\langle \partial_{\xi}^{-1} \partial_c u_{c(t)}, \tilde{u} \rangle_{L^2(\mathbb{R} \times \mathbb{T})} = 0
\right\}.
\end{equation}
Moreover, $c(t)$ and $\tilde{u}(t)$ satisfies
\begin{equation}
\label{orbit-bound-after}
|c(t) - c_*| +  \| \tilde{u}(t) \|_{H^1(\mathbb{R} \times \mathbb{T}) \cap H^1_{\mu}(\mathbb{R} \times \mathbb{T})}
\leq C \eps, \quad t \in \mathbb{R}_+.
\end{equation}
\end{lem}

\begin{Proof}
The proof is relatively well-known, see Proposition 5.1 in \cite{pw2}. It is based on the implicit function
theorem applied to the two constraints in the definition of $[X_{c(t)}^*]^{\perp}$ in (\ref{symplectic-orthogonality}).
Further details can be found in \cite{CCP}.
\end{Proof}

By the global well-posedness theory for the ZK equation (\ref{ZK}) \cite{Moliner},
there exists a unique global solution in class $u \in C(\mathbb{R}_+,H^1(\mathbb{R} \times \mathbb{T}) \cap H^1_{\mu}(\mathbb{R} \times \mathbb{T}))$
for $\mu > 0$ sufficiently small. By the initial bound (\ref{initial-bound}), the initial data
satisfy (\ref{orbit-bound}) for some $\eps > 0$ sufficiently small.
By the elementary continuation arguments, the decomposition (\ref{decomposition-lemma}) can be used
as long as the solution $u$ satisfies (\ref{orbit-bound}). The component $\tilde{u}$ in the decomposition
(\ref{decomposition-lemma}) satisfies the evolution equation (\ref{amplitude-evolution}) rewritten again as
\begin{equation}
\label{amplitude-evolution-0}
\tilde{u}_t = \partial_{\xi} (L_{c} - \partial_y^2 + 4(\dot{a}-c)) \tilde{u}
+ 4 (\dot{a} - c) \partial_{\xi} u_{c} - \dot{c} \partial_c u_{c}
- 6 \partial_{\xi} \tilde{u}^2,
\end{equation}
where the differential expression for $L_c$ is given by (\ref{Schrodinger-operator}).

Both parameters $a$ and $c$ depend on the time variable $t$.
Modulation equations for $a$ and $c$ are derived from the well-known projection algorithm, which has been
applied to similar problems in \cite{CCP,Marzuola,pw2}. The two constraints on $\tilde{u}$
in (\ref{symplectic-orthogonality}) represent the symplectic orthogonality conditions, which specify uniquely
$a$ and $c$ in the decomposition (\ref{decomposition-lemma}) as well as the time evolution of $a$ and $c$.
Moreover, one can show that $a,c \in C^1(\mathbb{R}_+)$.

If $\tilde{u}(0) \in [X_{c(0)}^*]^{\perp}$ initially, then $\tilde{u}(t)$ remains in $[X_{c(t)}^*]^{\perp}$
for every $t \in \mathbb{R}_+$, provided that the varying parameters $a,c \in C^1(\mathbb{R}_+)$
satisfy the system of modulation equations
\begin{eqnarray}
S \left[ \begin{array}{c} \dot{c} \\ 4 (\dot{a} - c) \end{array} \right]
= \frac{3}{\pi}  \left[ \begin{array}{c} \langle \partial_c u_c, \tilde{u}^2 \rangle_{L^2(\mathbb{R} \times \mathbb{T})} \\
\langle \partial_{\xi} u_c, \tilde{u}^2 \rangle_{L^2(\mathbb{R} \times \mathbb{T})}
\end{array} \right]\label{evolution-a-c}
\end{eqnarray}
with the coefficient matrix
\begin{equation}
\label{coefficeint-S}
S := \left[ \begin{array}{cc}
\frac{1}{2} (M'(c))^2 - \frac{1}{2\pi} \langle \partial_{\xi}^{-1} \partial^2_c u_c, \tilde{u} \rangle_{L^2(\mathbb{R} \times \mathbb{T})} &
P'(c) + \frac{1}{2\pi} \langle \partial_c u_c, \tilde{u} \rangle_{L^2(\mathbb{R} \times \mathbb{T})} \\
P'(c) - \frac{1}{2\pi} \langle \partial_c u_c, \tilde{u} \rangle_{L^2(\mathbb{R} \times \mathbb{T})} &
\frac{1}{2\pi} \langle \partial_{\xi} u_c, \tilde{u} \rangle_{L^2(\mathbb{R} \times \mathbb{T})} \end{array} \right],
\end{equation}
where $M(c) = \int_{\mathbb{R}} u_c(\xi) d\xi$ and $P(c) = \frac{1}{2} \int_{\mathbb{R}} u_c^2(\xi) d \xi$.
From the expression (\ref{line-soliton}), we obtain $M'(c) = 1/\sqrt{c}$ and $P'(c) = \sqrt{c}$.

\subsection{A secondary decomposition for $c = c_*$}

By Lemma \ref{theorem-bifurcation}, if $c = c_* = \frac{1}{5}$ and $\mu > 0$ is sufficiently small,
then $Y_{c_*} = {\rm span}\{ \psi_* \}$ is an invariant subspace of $L^2_{\mu}(\mathbb{R})$ for
the simple zero eigenvalue of the linearized operator
$\partial_{\xi} (L_{c_*} + 1) : H^3_{\mu}(\mathbb{R}) \to L^2_{\mu}(\mathbb{R})$.
Similarly,  $Y_{c_*}^* = {\rm span}\{ \eta_* \}$
is an invariant subspace of $L^2_{-\mu}(\mathbb{R})$ for the simple zero eigenvalue of the adjoint
operator $-(L_{c_*} + 1)\partial_{\xi} : H^3_{-\mu}(\mathbb{R}) \to L^2_{-\mu}(\mathbb{R})$.
We note the double degeneracy of the Fourier harmonics $e^{iy}$ and $e^{-iy}$,
when general transverse perturbations are considered.

The following lemma states the secondary decomposition of the solution $\tilde{u}$
defined in the primary decomposition (\ref{decomposition-lemma}).

\begin{lem}
\label{lem-decomposition-secondary}
Under assumptions of Lemma \ref{lem-decomposition}, let
$\tilde{u} \in C(\mathbb{R}_+,H^1(\mathbb{R} \times \mathbb{T}) \cap H^1_{\mu}(\mathbb{R} \times \mathbb{T}))$
be given by the decomposition (\ref{decomposition-lemma}) and (\ref{symplectic-orthogonality}).
There exist $b \in C(\mathbb{R}_+)$ and $v \in
C(\mathbb{R}_+,H^1(\mathbb{R} \times \mathbb{T}) \cap H^1_{\mu}(\mathbb{R} \times \mathbb{T}))$
such that the decomposition
\begin{equation}
\label{decomposition-secondary}
\tilde{u}(\xi,y,t) = \left( b(t) e^{iy} + \bar{b}(t) e^{-iy} \right) \psi_*(\xi) + v(\xi,y,t),
\end{equation}
holds with $v(t) \in [Y_{c(t)}^*]^{\perp}$ for every $t \in \mathbb{R}_+$, where
\begin{equation}
\label{symplectic-orthogonality-single}
[Y_{c(t)}^*]^{\perp} = \left\{ v \in [X_{c(t)}^*]^{\perp} : \quad
\langle \eta_* e^{iy}, v \rangle_{L^2(\mathbb{R} \times \mathbb{T})}
= \langle \eta_* e^{-iy}, v \rangle_{L^2(\mathbb{R} \times \mathbb{T})} = 0
\right\}.
\end{equation}
\end{lem}

\begin{Proof}
The proof is straightforward thanks to the fact
$\langle \eta_*, \psi_* \rangle_{L^2(\mathbb{R})} \neq 0$
by (\ref{limiting-coefficient}).
\end{Proof}

We introduce the decomposition
\begin{equation}
\label{decomposition-a-c}
a(t) = \int_0^t c(t') dt' + h(t), \quad c(t) = c_* + \delta(t),
\end{equation}
in addition to the decomposition (\ref{decomposition-secondary}).
To simplify notations, we also write $L_c = L_{c_*} + \Delta L_c$,
where $\Delta L_c \in L^{\infty}(\mathbb{R})$ satisfies the bound $\| \Delta L_c \|_{L^{\infty}} \leq A |c-c_*|$
for $|c-c_*|$ sufficiently small with a $c$-independent positive constant $A$.

The correction term $v$ in the decomposition (\ref{decomposition-secondary}) satisfies the time evolution equation
\begin{eqnarray}
\nonumber
v_t & = & \partial_{\xi} (L_{c_*} - \partial_y^2 + 4 \dot{h} + \Delta L_c) v
+ 4 \dot{h} \partial_{\xi} u_{c_*+\delta} - \dot{\delta} \partial_c u_{c_*+\delta} - (\dot{b}e^{iy} + \dot{\bar{b}} e^{-iy}) \psi_* \\
\label{amplitude-evolution-1}
& \phantom{t} &
+ \partial_{\xi} \left( 4 \dot{h} + \Delta L_c \right) (be^{iy} + \bar{b} e^{-iy}) \psi_*
- 6 \partial_{\xi} \left( (b e^{iy} + \bar{b} e^{-iy} ) \psi_* + v \right)^2.
\end{eqnarray}

The two constraints in (\ref{symplectic-orthogonality-single}) represent the symplectic orthogonality conditions,
which specify uniquely the complex parameter $b$ in the secondary decomposition (\ref{decomposition-secondary}).
Again, one can show that $b \in C^1(\mathbb{R}_+)$

If $v(0) \in [Y_{c(0)}^*]^{\perp}$ initially, then $v(t)$ remains in $[Y_{c(t)}^*]^{\perp}$ for every $t \in \mathbb{R}_+$,
provided that the varying parameter $b \in C^1(\mathbb{R}_+)$ satisfy the following modulation equation:
\begin{eqnarray}
\nonumber
& \phantom{t} & \dot{b} \langle \eta_*, \psi_* \rangle_{L^2(\mathbb{R})} +
b \langle \psi_*, (4 \dot{h} + \Delta L_c) \psi_* \rangle_{L^2(\mathbb{R})} +
\frac{1}{2\pi} \langle \psi_* e^{iy}, \Delta L_c v \rangle_{L^2(\mathbb{R} \times \mathbb{T})}\\
& \phantom{t} & \phantom{texttexttext}
= \frac{3}{\pi}  \langle \psi_* e^{iy}, \left[ (b e^{iy} + \bar{b} e^{-iy} ) \psi_* + v \right]^2 \rangle_{L^2(\mathbb{R} \times \mathbb{T})}.
\label{evolution-b}
\end{eqnarray}
Substituting (\ref{decomposition-secondary}) and (\ref{decomposition-a-c})
to the system (\ref{evolution-a-c}) yields the equivalent form of the modulation equations:
\begin{eqnarray}
S \left[ \begin{array}{c} \dot{\delta} \\ 4 \dot{h} \end{array} \right]
= \frac{3}{\pi}
\left[ \begin{array}{c} \langle \partial_c u_c, \left[ (b e^{iy} + \bar{b} e^{-iy} ) \psi_* + v \right]^2 \rangle_{L^2(\mathbb{R} \times \mathbb{T})} \\
\langle \partial_{\xi} u_c, \left[ (b e^{iy} + \bar{b} e^{-iy} ) \psi_* + v \right]^2 \rangle_{L^2(\mathbb{R} \times \mathbb{T})}
\end{array} \right],\label{evolution-ac}
\end{eqnarray}
where $c(t) = c_* + \delta(t)$ and $S$ in (\ref{coefficeint-S}) becomes now
\begin{equation}
\label{coefficeint-S-S}
S := \left[ \begin{array}{cc}
\frac{1}{2} (M'(c))^2 - \frac{1}{2\pi} \langle \partial_{\xi}^{-1} \partial^2_c u_c, v \rangle_{L^2(\mathbb{R} \times \mathbb{T})} &
P'(c) + \frac{1}{2\pi} \langle \partial_c u_c, v \rangle_{L^2(\mathbb{R} \times \mathbb{T})} \\
P'(c) - \frac{1}{2\pi} \langle \partial_c u_c, v \rangle_{L^2(\mathbb{R} \times \mathbb{T})} &
\frac{1}{2\pi} \langle \partial_{\xi} u_c, v \rangle_{L^2(\mathbb{R} \times \mathbb{T})} \end{array} \right].
\end{equation}
The system (\ref{evolution-b}) and (\ref{evolution-ac}) determine the time evolution of the varying
parameters $b$, $h$, and $\delta$, whereas the evolution
problem (\ref{amplitude-evolution-1}) determines the correction term $v(t) \in  [Y_{c(t)}^*]^{\perp}$.

\subsection{Justification of the approximation error}

We justify the error bound (\ref{final-bound}) with the energy method pioneered in \cite{pw2}.
First, we recall that the energy
\begin{equation}
\label{energy}
E(u) = \frac{1}{2} \int_{\mathbb{R} \times \mathbb{T}} \left[ u_x^2 + u_y^2 - 4 u^3 \right] dx dy
\end{equation}
and the momentum
\begin{equation}
\label{momentum}
Q(u) = \frac{1}{2} \int_{\mathbb{R} \times \mathbb{T}} u^2 dx dy
\end{equation}
are conserved in time $t$ for a global solution $u \in C(\mathbb{R},H^1(\mathbb{R}\times\mathbb{T}))$ to the ZK equation (\ref{ZK}).
The line soliton (\ref{line-soliton}) is a critical point of the action functional
$\Lambda_c(u) := E(u) + 4c Q(u)$, see (\ref{action}). Thanks to the translational invariance
of the ZK equation (\ref{ZK}), the decomposition (\ref{decomposition-lemma}) yields
\begin{equation}
\label{decomposition-Lambda}
\Lambda_c(u_c + \tilde{u}) - \Lambda_c(u_c) = \frac{1}{2} \langle (L_c - \partial_y^2) \tilde{u}, \tilde{u} \rangle_{L^2}
+ N_c(\tilde{u}),
\end{equation}
where the differential expression for $L_c$ is given by (\ref{Schrodinger-operator}) and $N_c$ is a nonlinear term satisfying
\begin{equation}
\label{N-c}
\left| N_c(\tilde{u}) \right| \leq A \| \tilde{u} \|_{H^1}^3,
\end{equation}
for some positive constant $A$ as long as $\| \tilde{u} \|_{H^1}$ is small.
By using the Fourier series
$$
\tilde{u}(\xi,y,t) = \frac{1}{\sqrt{2\pi}} \sum_{k \in \mathbb{Z}} \hat{u}_k(\xi,t) e^{iky}
$$
and Parseval's equality, we can represent the second variation of $\Lambda_c$ at $u_c$ in the form
\begin{equation}
\label{second-variation}
\langle (L_c - \partial_y^2) \tilde{u}, \tilde{u} \rangle_{L^2(\mathbb{R} \times \mathbb{T})} =
\sum_{k\in \mathbb{Z}} \langle (L_c + k^2) \hat{u}_k, \hat{u}_k \rangle_{L^2(\mathbb{R})}.
\end{equation}
The following lemma summarizes the coercivity results for the second variation of $\Lambda_c$ at $u_c$.

\begin{lem}
For $\mu > 0$ sufficiently small, there exists a constant $A > 0$ such that
for every $\hat{u}_k \in H^1(\mathbb{R})$ and for every $k \in \mathbb{Z} \backslash \{0,\pm 1\}$,
it is true that
\begin{equation}
\label{coercivity-3}
\langle (L_{c_*} + k^2)\hat{u}_k,\hat{u}_k \rangle_{L^2(\mathbb{R})} \geq A \| \hat{u}_k \|_{H^1(\mathbb{R})}^2,
\end{equation}
whereas for every $\hat{u}_0, \hat{u}_{\pm 1} \in H^1(\mathbb{R}) \cap H^1_{\mu}(\mathbb{R})$, it is true that
\begin{equation}
\label{coercivity-1}
\langle L_{c_*} \hat{u}_0,\hat{u}_0 \rangle_{L^2(\mathbb{R})} \geq A \| \hat{u}_0 \|_{H^1(\mathbb{R})}^2
\quad \mbox{\rm if} \;\; \langle u_{c_*}, \hat{u}_0 \rangle_{L^2(\mathbb{R})} =
\langle \partial_{\xi}^{-1} \partial_c u_{c} |_{c = c_*}, \hat{u}_0 \rangle_{L^2(\mathbb{R})} = 0
\end{equation}
and
\begin{equation}
\label{coercivity-2}
\langle (L_{c_*} + 1) \hat{u}_{\pm 1},\hat{u}_{\pm 1} \rangle_{L^2(\mathbb{R})} \geq A \| \hat{u}_{\pm 1} \|_{H^1(\mathbb{R})}^2
\quad \mbox{\rm if} \;\; \langle \eta_*, \hat{u}_{\pm 1} \rangle_{L^2(\mathbb{R})} = 0.
\end{equation}
\label{lemma-coercivity}
\end{lem}

\begin{Proof}
The spectral information on the Schr\"{o}dinger operator $L_c : H^2(\mathbb{R}) \to L^2(\mathbb{R})$
with the two lowest eigenvalues (\ref{first-eig}) and (\ref{second-eig}) is sufficient to conclude that
$L_c + k^2 : H^2(\mathbb{R}) \to L^2(\mathbb{R})$ is strictly positive for every $k \in \mathbb{Z} \backslash \{0,\pm 1\}$.
The bound (\ref{coercivity-3}) follows by the spectral theorem and G{\aa}rding's inequality.

Since $L_{c_*} + 1 : H^2(\mathbb{R}) \to L^2(\mathbb{R})$ is non-negative with a one-dimensional kernel spanned by
$\psi_*$ and $\langle \eta_*, \psi_* \rangle_{L^2(\mathbb{R})} \neq 0$ by (\ref{limiting-coefficient}),
this operator is strictly positive under the constraint in (\ref{coercivity-2}). The constraint
in (\ref{coercivity-2}) is well-defined if $\hat{u}_{\pm 1} \in H^1_{\mu}(\mathbb{R})$.
The bound (\ref{coercivity-2}) follows by the spectral theorem and G{\aa}rding's inequality.

Since $L_{c_*} : H^2(\mathbb{R}) \to L^2(\mathbb{R})$ has one negative and one simple eigenvalues, whereas $P'(c_*) > 0$,
this operator is non-negative under the first constraint in (\ref{coercivity-1})
with a one-dimensional kernel spanned by $\partial_{\xi} u_{c_*}$ \cite{pw1}.
Since
$$
\langle u_{c_*}, \partial_{\xi} u_{c_*} \rangle_{L^2(\mathbb{R})} = 0, \quad
\langle \partial_{\xi}^{-1} \partial_c u_{c} |_{c = c_*}, \partial_{\xi} u_{c_*} \rangle_{L^2(\mathbb{R})} \neq 0,
$$
this operator is strictly positive under the two constraints in (\ref{coercivity-1}).
The second constraint
in (\ref{coercivity-1}) is well-defined if $\hat{u}_0 \in H^1_{\mu}(\mathbb{R})$.
The bound (\ref{coercivity-1}) follows by the spectral theorem and G{\aa}rding's inequality.
\end{Proof}

In order to justify the error bound (\ref{final-bound}), we construct the following
energy function
\begin{eqnarray}
\label{energy-function-F}
F(c) :=  E(u) - E(u_{c_*}) + 4 c \left[ Q(u) - Q(u_{c_*}) \right].
\end{eqnarray}
Thanks to the conservation of energy $E$ and momentum $Q$ in time $t$,
we have
\begin{eqnarray}
\label{energy-function-F-time-0}
F(c) = E(u_0) - E(u_{c_*}) + 4c \left[ Q(u_0) - Q(u_{c_*}) \right]
\end{eqnarray}
Since $c(t)$ depends on $t$, $F(c(t))$ depends on $t$ but only linearly in $c(t)$.
By Lemmas \ref{lem-decomposition} and \ref{lem-decomposition-secondary}, we rewrite
the decompositions (\ref{decomposition-lemma}), (\ref{decomposition-secondary}), and (\ref{decomposition-a-c}) in the form
\begin{equation}
\label{decomposition-bound}
u(x,y,t) = u_{c_* + \delta(t)}(\xi) + \left( b(t) e^{iy} + \bar{b}(t) e^{-iy} \right) \psi_*(\xi) + v(\xi,y,t).
\end{equation}
Substituting the decomposition (\ref{decomposition-bound}) into (\ref{energy-function-F}) yields the following
\begin{equation}
\label{decomposition-F}
F(c) = D(c) + \frac{1}{2} \langle (L_c - \partial_y^2) v,v \rangle_{L^2} + N_c((b e^{iy} + \bar{b} e^{-iy}) \psi_* + v),
\end{equation}
where the expansion (\ref{decomposition-Lambda}) has been used, $N_c$ satisfies (\ref{N-c}), and
\begin{eqnarray}
\nonumber
D(c) & := & E(u_c) - E(u_{c_*}) + 4 c \left[ Q(u_c) - Q(u_{c_*}) \right] \\
& = & \frac{1}{2} D''(c_*) (c-c_*)^2 + \tilde{D}(c).\label{decomposition-D}
\end{eqnarray}
The latter expansion is obtained from $D(c_*) = D'(c_*) = 0$ and $D''(c_*) = 4 P'(c_*) > 0$,
where $P(c) = \frac{1}{2} \int_{\mathbb{R}} u_c^2(\xi) d \xi$, thanks
to the variational characterization of the line soliton (\ref{line-soliton}) with the action
functional (\ref{action}). Thanks to the smoothness of $D$ in $c$, we have
$\tilde{D}(c) = \mathcal{O}((c-c_*)^3)$ as $c \to c_*$.

The following result transfers smallness of the initial bound (\ref{initial-bound}) to
smallness of $F(c_*)$ and $Q(u_0) - Q(u_{c_*})$ in (\ref{energy-function-F-time-0}).

\begin{lem}
\label{lem-energy}
There exists an $\eps$-independent positive constant $A$ such that
\begin{equation}
\label{bound-F-P}
|F(c_*)| \leq A \eps^4, \quad |Q(u_0) - Q(u_{c_*})| \leq A  \eps^2.
\end{equation}
\end{lem}

\begin{Proof}
The second bound in (\ref{bound-F-P}) follows from the initial bound (\ref{initial-bound})
thanks to the bound $|c(0) - c_*| \leq \eps^2$ and the triangle inequality.
The first bound in (\ref{bound-F-P}) follows from the expansion (\ref{decomposition-F}) with $D(c_*) = 0$,
the cubic term vanishing
\begin{equation}
\label{cubic-term}
N_c((b e^{iy} + \bar{b} e^{-iy}) \psi_*) = 0,
\end{equation}
the definition $b(0) = \eps$, and the triangle inequality.
\end{Proof}

Under the two constraints in (\ref{symplectic-orthogonality}) and the two constraints in (\ref{symplectic-orthogonality-single}),
it follows from Lemma \ref{lemma-coercivity} that there exists an $\eps$-independent constant $A$ such that
\begin{equation}
\label{lower-bound}
\langle (L_c - \partial_y^2) v, v \rangle_{L^2} \geq A \| v \|^2_{H^1(\mathbb{R} \times \mathbb{T})}.
\end{equation}

Let us assume that
\begin{equation}
\label{assumption-on-b}
|b(t)| \leq C_0 \eps, \quad t \in \mathbb{R}_+.
\end{equation}
This assumption is true at $t = 0$ since $b(0) = \eps$ and it remains true
for every $t \in \mathbb{R}_+$ as long as the normal form (\ref{normal-form})
with $|c_+ - c_* | \leq C_0 \eps^2$ can be used.

Combining (\ref{energy-function-F-time-0}), (\ref{decomposition-F}), (\ref{decomposition-D}),
and (\ref{lower-bound}) yields the following lower bound:
\begin{eqnarray*}
F(c_*) & \geq & 2 P'(c_*) (c-c_*)^2 - 4 (c-c_*) [ Q(u_0) - Q(u_{c_*})] + \tilde{D}(c) \\
& \phantom{t} &
+ \frac{1}{2} A \| v \|^2_{H^1(\mathbb{R} \times \mathbb{T})} + N_c((b e^{iy} + \bar{b} e^{-iy}) \psi_* + v).
\end{eqnarray*}
Thanks to the bounds (\ref{bound-F-P}) in Lemma \ref{lem-energy}, the smallness of $N_c$ in (\ref{N-c}),
and the cubic term vanishing in (\ref{cubic-term}), we obtain the bound (\ref{final-bound}) for $|c(t)-c_*|$
and $\| v \|_{H^1(\mathbb{R} \times \mathbb{T})}$.

The bound (\ref{final-bound}) on $|\dot{a}(t) - c(t)|$ follows by the expansion (\ref{decomposition-a-c})
as long as the assumption (\ref{assumption-on-b}) is true, since $\dot{h} = \mathcal{O}(|b|^2)$
follows from the modulation equations (\ref{evolution-ac}), see estimate (\ref{estimates-on-b-h}) below.

Thus, Theorem \ref{theorem-normal-form} is proven as long as the normal form (\ref{normal-form})
is derived and justified. This will be done with the near-identity transformations and the momentum conservation.

\subsection{Near-identity transformations}

Because $u_{c_*}$ and $\psi_*^2$ are even functions of $\xi$, whereas $P'(c_*) \neq 0$,
the modulation equations (\ref{evolution-b}) and (\ref{evolution-ac})
yields the following balance at the leading order:
\begin{equation}
\label{estimates-on-b-h}
\dot{b} = \mathcal{O}((|\delta| + |b|^2) |b|), \quad \dot{h} = \mathcal{O}(|b|^2),
\end{equation}
whereas the source terms in the evolution problem (\ref{amplitude-evolution-1}) are of the order
of $\mathcal{O}(|b|^2)$. In what follows, we write out the leading-order terms provided that $\delta$ and $b$
remain small for every $t \in \mathbb{R}_+$. Recall that the initial bound (\ref{initial-bound})
yields $|\delta(0)| \leq \eps^2$ and $b(0) = \eps$, where $\eps \in (0,\eps_0)$ is a small parameter.
Smallness of $\delta(t) := c(t) - c_*$ for every $t \in \mathbb{R}_+$ is proven in Section 3.3,
see the bound (\ref{final-bound}). Smallness of $b(t)$ for every $t \in \mathbb{R}_+$
is assumed in the bound (\ref{assumption-on-b}) and is proven here from the normal form (\ref{normal-form}).

In order to derive the normal form (\ref{normal-form}), we use the near-identity transformations, which
are very similar to the ones used in the proof of Lemma \ref{lemma-modulated-wave}.
In particular, we will remove the $\mathcal{O}(|b|^2)$
terms in the equation for $\dot{h}$ and $v_t$. Hence, we represent
the correction term $v$ in the decomposition (\ref{decomposition-secondary}) as follows:
\begin{equation}
\label{decomposition-tertiary}
v(\xi,y,t) = \left( b(t)^2 e^{2iy} + \bar{b}(t)^2 e^{-2iy} \right) w_2(\xi) +
|b(t)|^2 w_0(\xi) + w(\xi,y,t),
\end{equation}
where $w_0$ and $w_2$ are the same solutions of the linear inhomogeneous equations (\ref{inhom-eq-1})
and (\ref{inhom-eq-2}), whereas $w$ satisfies the transformed evolution equation
\begin{eqnarray}
\label{amplitude-evolution-2}
w_t & = & \partial_{\xi} (L_{c_*} - \partial_y^2 + 4 \dot{h} + \Delta L_c) w
+ 4 \dot{h} \partial_{\xi} u_{c_*+\delta} - \dot{\delta} \partial_c u_{c_*+\delta}  \\
\nonumber
& \phantom{t} &
- (\dot{b}e^{iy} + \dot{\bar{b}} e^{-iy}) \psi_*
- (2 b \dot{b} e^{2iy} + 2 \bar{b} \dot{\bar{b}} e^{-2iy}) w_2 - (\bar{b} \dot{b} + b \dot{\bar{b}}) w_0 \\
\nonumber
& \phantom{t} &
+ \partial_{\xi} \left( 4 \dot{h} + \Delta L_c \right)
\left[ (be^{iy} + \bar{b} e^{-iy}) \psi_* + ( b^2 e^{2iy} + \bar{b}^2 e^{-2iy}) w_2 +
|b|^2 w_0 \right] \\
\nonumber
& \phantom{t} &
- 12 \partial_{\xi} (b e^{iy} + \bar{b} e^{-iy} ) \psi_* \left[( b^2 e^{2iy} + \bar{b}^2 e^{-2iy}) w_2 +
|b|^2 w_0\right] \\
\nonumber
& \phantom{t} &
- 6 \partial_{\xi} \left[(b^2 e^{2iy} + \bar{b}^2 e^{-2iy}) w_2 +
|b|^2 w_0 \right]^2.
 \end{eqnarray}
We rewrite the first equation in the system (\ref{evolution-ac}) as follows:
\begin{equation}
\label{normal-form-1a}
4 P'(c_*) \dot{h} = 12 |b|^2 \langle \partial_c u_{c_*}, \psi_*^2 \rangle_{L^2} + \mathcal{O}(|b|^4)
= \frac{48}{5 \sqrt{c_*}} |b|^2 + \mathcal{O}(|b|^4),
\end{equation}
where the explicit expression  (\ref{explicit-solution-2}) has been used.
Since $P'(c_*) = \sqrt{c_*}$ and $c_* = \frac{1}{5}$, we obtain
\begin{equation}
\label{decomposition-h-dot}
\dot{h} = 12 |b|^2 + \mathcal{O}(|b|^4)
\end{equation}
and
\begin{equation}
\label{decomposition-last}
w(\xi,y,t) = 12 |b|^2 \partial_{c} u_{c_*}(\xi) + \tilde{w}(\xi,y,t),
\end{equation}
where $\tilde{w}$ satisfied a transformed evolution equation without the $\mathcal{O}(|b|^2)$ terms
in the right-hand side of (\ref{amplitude-evolution-2}).
Substituting (\ref{decomposition-tertiary}), (\ref{decomposition-h-dot}), and (\ref{decomposition-last})
into the modulation equation (\ref{evolution-b}) yields
\begin{eqnarray}
\label{normal-form-1}
\dot{b} \langle \eta_*, \psi_* \rangle_{L^2} & = &
12 |b|^2 b \langle \psi_*^2, w_0 + w_2 \rangle_{L^2} +
144 |b|^2 b \langle \psi_*^2, \partial_c u_{c_*} \rangle_{L^2} \\
\nonumber
& \phantom{t} & - b \langle \psi_*, (L_{c_*}' \delta + 48 |b|^2) \psi_* \rangle_{L^2}
+ \mathcal{O}(\delta^2 |b| + |b|^5),
\end{eqnarray}
where $L_{c_*}'$ is given by (\ref{derivative-operator}) and
we have used $\Delta L_c = L_{c_*}' \delta + \mathcal{O}(\delta^2)$
for $\delta = c - c_*$. After straightforward computations,
equation (\ref{normal-form-1}) takes the form
\begin{equation}
\label{normal-form-1-equiv}
\dot{b} \langle \eta_*, \psi_* \rangle_{L^2}  =
-12 |b|^2 b \langle \psi_*^2, \tilde{w}_2 \rangle_{L^2} +
\frac{64}{\sqrt{c_*}} |b|^2 b + \frac{16}{3 \sqrt{c_*}} b \delta + \mathcal{O}(\delta^2 |b| + |b|^5),
\end{equation}
where $\tilde{w}_2$ is found from the solution of the linear inhomogeneous equation
(\ref{inhomogeneous-eq-tilde-w-2}).

The modulation equation (\ref{normal-form-1-equiv}) is not closed on $b$
because $\delta$ is related to $|b|^2$ by the second equation
of the system (\ref{evolution-ac}). In fact, this equation relates $\dot{\delta}$
to $\mathcal{O}(|\dot{\bar{b}} b|) = \mathcal{O}(|\delta| |b|^2 + |b|^4)$, however,
it yields $\delta = \mathcal{O}(|b|^2)$ after integration. In order to avoid
integration of the second equation of the system (\ref{evolution-a-c}), we
use the momentum conservation $Q(u) = Q(u_0)$, where the momentum $Q$ is given by
(\ref{momentum}).

Substituting decompositions (\ref{decomposition-lemma}), (\ref{decomposition-secondary}),
(\ref{decomposition-a-c}), (\ref{decomposition-tertiary}), and (\ref{decomposition-last}) into (\ref{momentum}) yields
the expansion
\begin{equation}
\label{expansion-Q}
Q(u) = 2 \pi \left[ P(c_* + \delta) + |b|^2 \| \psi_* \|^2_{L^2} + \mathcal{O}(|\delta| |b|^2 + |b|^4) \right],
\end{equation}
where $P(c) = \frac{1}{2} \int_{\mathbb{R}} u_c^2(\xi) d \xi$ and we have used the fact
$$
\langle u_{c_*}, w_0 \rangle_{L^2} + 12 \langle u_{c_*}, \partial_{c} u_{c_*} \rangle_{L^2} = 0,
$$
which follows from integration of the explicit expressions (\ref{line-soliton}), (\ref{explicit-solution-2}),
and (\ref{explicit-solution-1}). By computing $P(c)$ and $\| \psi_* \|^2_{L^2}$
from (\ref{line-soliton}) and (\ref{limiting-function}), we use the momentum conservation
and expand (\ref{expansion-Q}) to the explicit form:
\begin{equation}
\label{expansion-Q-Q}
Q(u_0) = 2 \pi \left[ P(c_*) + \sqrt{c_*} \delta + \frac{16}{15 \sqrt{c_*}} |b|^2
+ \mathcal{O}(\delta^2 + |b|^4) \right],
\end{equation}
which yields with $c_* = \frac{1}{5}$,
\begin{equation}
\label{normal-form-1e}
\delta = \delta_0 - \frac{16}{3} |b|^2 + \mathcal{O}(\delta_0^2 + |b|^4),
\end{equation}
where $\delta_0$ is a constant in $t$ determined by the initial data.

Substituting equation (\ref{normal-form-1e})
into the modulation equation (\ref{normal-form-1-equiv}) yield
\begin{equation}
\label{normal-form-2}
\dot{b} \langle \eta_*, \psi_* \rangle_{L^2}  =
-12 |b|^2 b \langle \psi_*^2, \tilde{w}_2 \rangle_{L^2} +
\frac{16}{3 \sqrt{c_*}} \left( \delta_0 + \frac{20}{3} |b|^2 \right) b +
\mathcal{O}(\delta_0^2 |b| + |b|^5).
\end{equation}
Defining $\delta_0 := c_+ - c_*$, using the explicit expression
(\ref{limiting-coefficient}) and (\ref{derivative-eigenvalue}),
and truncating (\ref{normal-form-2}) yield the normal form (\ref{normal-form}) with
\begin{equation}
\label{gamma}
\gamma := \frac{96}{\pi^2 \sqrt{c_*}} \left( - \frac{1}{5} \sqrt{c_*} \langle \psi_*^2, \tilde{w}_2 \rangle_{L^2} + \frac{16}{27} \right).
\end{equation}
Solutions to the normal form (\ref{normal-form}) under the assumption that $b(0) = \eps$ and
$|\delta_0| \leq C_0 \eps^2$ for an $\eps$-independent positive constant $C_0$ satisfy
the bound $|b(t)| \leq C_0 \eps$ for all $t \in \mathbb{R}_+$,
which is the same as the one used in (\ref{assumption-on-b}). For such solutions, the remainder term in (\ref{normal-form-2})
is of the order $\mathcal{O}(\eps^5)$, hence the truncation of (\ref{normal-form-2}) into (\ref{normal-form}) is justified
within the approximation error controlled by the bound (\ref{final-bound}).
An elementary continuation argument completes the proof of Theorem \ref{theorem-normal-form}.

\begin{rem}
\label{remark-coefficients}
The numerical coefficients $(\lambda'(c_*),\gamma)$ are different from the numerical
coefficients $(\alpha,\beta)$. This difference is explained
as follows. Expansions (\ref{decomposition-a-c}), (\ref{decomposition-h-dot}), and (\ref{normal-form-1e}) yield
$$
c = c_* + \delta_0 - \frac{16}{3} |b|^2 + \mathcal{O}(|b|^4), \quad
\dot{a} = c_* + \delta_0 + \frac{20}{3} |b|^2 + \mathcal{O}(|b|^4).
$$
The normal form (\ref{normal-form-2}) in the stationary case $\dot{b} = 0$
corresponds to the effective speed correction given by
$$
\dot{a} - c_* =  \delta_0 + \frac{20}{3} |b|^2 + \mathcal{O}(|b|^4) =
\frac{9}{4} \sqrt{c_*} |b|^2 \langle \psi_*^2, \tilde{w}_2 \rangle_{L^2} + \mathcal{O}(|b|^4)
= -\frac{\beta}{\alpha} |b|^2 + \mathcal{O}(|b|^4),
$$
in agreement with the cubic algebraic equation (\ref{normal-form-static}).
In the time-dependent case, the roles of $\delta_0$ and $\frac{20}{3} |b|^2$ are different
because the former is constant in $t$ but the latter changes in $t$. A very similar discrepancy
between numerical coefficients of the stationary and time-independent normal forms
is observed in \cite{Phan} in the context of symmetry-breaking bifurcations in a double-well
potential.
\end{rem}

\begin{rem}
\label{remark-numerics}
Bifurcation analysis of Lemma \ref{lemma-modulated-wave} relies on the fact that
the coefficient $\beta$ of the cubic term in the normal form (\ref{normal-form-static})
is negative. This fact has been proven in (\ref{gamma-stationary}).
It is equally important for the stability analysis near the line soliton
that the coefficient $\gamma$ of the cubic term in the normal form (\ref{normal-form}) is negative.
Since $\tilde{w}_2(\xi) \geq 0$ for all $\xi \in \mathbb{R}$, the first term in $\gamma$ is negative.
On the other hand, the second term is positive, so that
$\gamma < 0$ only if the negative term prevails. We show this fact with
the following numerical computation.

We approximate the function $\tilde{w}_2$ by using the central-difference method for
the linear inhomogeneous equation (\ref{inhomogeneous-eq-tilde-w-2})
and then approximate the integral $\langle \psi_*^2, \tilde{w}_2 \rangle_{L^2}$
by using the composite trapezoidal method.
Testing the codes on the function $w_0$ which satisfies the linear inhomogeneous equation (\ref{inhom-eq-1})
with the explicit solution (\ref{explicit-solution-1}) yields
$$
\frac{1}{5} \sqrt{c_*} \langle \psi_*^2, w_0 \rangle_{L^2} \approx -1.5238
$$
which corresponds to the exact value $-\frac{32}{21}$ within the computational error of $\mathcal{O}(10^{-5})$. Performing the same
task for $\tilde{w}_2$ which satisfies the linear inhomogeneous equation (\ref{inhomogeneous-eq-tilde-w-2}) yields
$$
\frac{1}{5} \sqrt{c_*} \langle \psi_*^2, \tilde{w}_2 \rangle_{L^2} \approx 1.2359,
$$
which is essentially bigger than $\frac{16}{27} \approx 0.5926$. Therefore,
$$
\frac{1}{5} \sqrt{c_*} \langle \psi_*^2, \tilde{w}_2 \rangle_{L^2} > \frac{16}{27},
$$
which implies that $\gamma < 0$ in (\ref{gamma}).
\end{rem}

\end{document}